\documentclass[10pt]{article}

\usepackage[dvips]{graphicx}

\begin{document}

\title{Simulating Protein Conformations through \\ Global Optimization}

\author{A. Mucherino\thanks{Laboratorie d'Informatique (LIX), \'Ecole Polytechnique, Palaiseau, France} \and
        O. Seref\thanks{Department of Business Information Technology, Virginia Polytechnic Institute and State University, Blacksburg, VA, USA} \and
        P.M. Pardalos\thanks{Center for Applied Optimization, University of Florida, FL, USA}
}

\maketitle

%%%

\begin{abstract}
Many researches have been working on the protein folding problem from more than half century. Protein folding is indeed one of the major unsolved problems in science. In this work, we discuss a model for the simulation of protein conformations. This simple model is based on the idea of imposing few geometric requirements on chains of atoms representing the backbone of a protein conformation. The model leads to the formulation of a global optimization problem, whose solutions correspond to conformations satisfying the desired requirements. The global optimization problem is solved by the recently proposed Monkey Search algorithm. The simplicity of the optimization problem and the effectiveness of the used meta-heuristic search allowed the simulation of a large set of high-quality conformations. We show that, even though only few geometric requirements are imposed, some of the simulated conformation results to be similar (in terms of RMSD) to conformations real proteins actually have in nature.
%\keywords{protein simulations \and protein folding \and global optimization \and meta-heuristics}
\end{abstract}

%%%

%%%%%%%%%%%%%%%%%%%%%%%%%%%%%%%%%%%%%%%%%%%%%%%%%%%%%%%%%%%%%%%%%%%%%%%%%%%%%%%%%%%
\section{Introduction}  \label{Introduction}

Protein folding is one of the major challenges in science \cite{Science}. Living organisms contain plenty of proteins, which perform many functions of vital importance. Proteins are chains of smaller molecules called amino acids. They may be formed by either one chain only, or more chains. The sequence of 20 amino acids along the protein chains is considered to contain all the information needed for the protein to fold into its functional structure \cite{Anfinsen}. The importance the scientific community is giving to proteins can be simply evaluated by performing a search of papers with Google Scholar: more than 5000 papers published between January and October 2008 contain the words ``protein folding''.

The conformation which a protein have in a cell defines its dynamics in the cell and therefore its function. Experimental techniques, such as X-ray crystallography and Nuclear Magnetic Resonance (NMR), are able to provide information about the protein conformations. Unfortunately, these techniques can be quite expensive and require a lot of time. The two techniques have different advantages and disadvantages, and they are usually complementary \cite{Yee}. Therefore, many researchers have been working with the aim of understanding the folding process and reproducing it artificially on computers. This is one of the subproblems in which the protein folding problem can be divided \cite{Dill}. Many algorithms and computational procedures have been developed over the years for attempting the prediction of the protein conformations. Some methods are based on the knowledge derived from protein databases, where as some are solely based on the idea that the energy involved in the protein molecule must be as minimum as possible in the stable and functional protein conformations. Different force fields have been proposed for modeling the energy inside a protein. Examples of force fields are AMBER \cite{Duan}, CHARMM22 \cite{MacKerell}, GROMOS \cite{vanGunsteren}, ECEPP/3 \cite{Nemethy} and ECEPP/5 \cite{Arnautova}. These force fields differ from each other regarding the kind of interactions taken into account and for the representation used for modeling the protein molecules.

The unbounded interactions, which involve atoms that are not chemically bonded and that may be far on the protein chain, are among the interactions usually considered. Such interactions involve the electrostatic forces and the so-called Van der Waals forces. The first ones simply model the attractive forces between atoms having different charges and the repulsive forces between two atoms with the same charge. The Van der Waals forces are instead more complex, and they are usually modeled by the Lennard Jones energy function \cite{LennardJones}. All of these interactions can be considered at an atomic level. However, this approach is complex and difficult to model, especially for large proteins with thousands of atoms.

The Lennard Jones energy is also used for finding optimal clusters of molecules. This is a well-known difficult problem in global optimization due to the large number of local minima \cite{Northby,Xiang}. The problem of predicting the conformation of a protein is therefore much more complex, since the Lennard Jones energy is usually only a term of the objective function to be minimized when dealing with protein predictions. Despite all efforts in this field, the research community is still far from finding a good solution to the problem of protein prediction with purely physics-based methods.

Within the last decade there has been a new emerging approach to protein folding. Baker in \cite{Baker} suggested that the protein structures are dependent on the topology of their native states and that the specific details of the interactions stabilizing these structures play a minor role in the folding process. This new emerging approach makes a large use of geometric properties and features of the protein conformations. The studies presented in \cite{Banavar04,Banavar03,Banavar02,Hoang} are first examples, in which this new approach has been used. In this work, we mainly refer to the studies described in \cite{Ceci,MucherinoPF}, where protein conformations have been simulated by using a model mainly based on geometric properties. 

This model represents a different way to consider the problem. The aim is not to predict the protein conformation related to a specific protein sequence, but rather to simulate protein conformations having some predetermined properties. In this paper, we present a very simple model based on geometric properties only. We impose few geometric requirements on a chain of amino acids and we simulate conformations satisfying such requirements. Each geometric requirement is controlled through a simple function: the more this function is minimized, the more the amino acids satisfy the corresponding requirements. Our model leads to a global optimization problem, in which the objective function is composed of the weighted sum of as many terms as the number of imposed requirements. 

We simulate a large subset of conformations and we compare these conformations to the ones of target proteins. In order to generate this large subset in a reasonable amount of time, we employ the recently proposed Monkey Search \cite{MucherinoMS} to solve the global optimization problem. This meta-heuristic method has been shown to be more efficient compared to the other meta-heuristic search methods on difficult optimization problems. We show that some simulated conformations are similar to real protein conformations chosen as targets. The comparisons are made in terms of Root Mean Square Deviation (RMSD) values.

The rest of the paper is organized as follows. In Section \ref{sec:opt} we present our model for protein simulations in detail, and formulate a global optimization problem based on this model. In Section \ref{sec:sim} we  show the differences and similarities between our model and other models in the literature. In Section \ref{sec:ms} the meta-heuristic method used for solving the optimization problem is discussed. In Section \ref{sec:experim} computational experiments are shown and conclusions are drawn in Section \ref{sec:concl}.

%%%%%%%%%%%%%%%%%%%%%%%%%%%%%%%%%%%%%%%%%%%%%%%%%%%%%%%%%%%%%%%%%%%%%%%%%%%%%%%%%%%
\section{A simple geometric model}  \label{sec:opt}

In this section, we present the details of a simple model for the simulation of protein conformations. The model considers geometric properties only. The protein conformations $X = \{ x_1,x_2,\dots,x_n \}$ are represented through the spatial positions $x_i$ of their $C_{\alpha}$ Carbon atoms. Information regarding the side chains of amino acids are not considered at all, and the amino acid is represented by its $C_{\alpha}$ Carbon atom only. We impose few geometric requirements on this simple representation and investigate how they can influence the simulation of conformations with the typical properties of the real conformations. 

The first requirement we introduce is in regards to the compactness of the protein conformations. Globular proteins are compact in general. Water, which is a protein solvent, is not allowed among the amino acids, except for small cavities in proteins that are filled with few water molecules. Therefore, all of the amino acids in a protein molecule must be as close to each other as possible. This means that the Euclidean distances between the amino acids should be minimized. Note that the consecutive amino acids have a constant distance because of the peptide bond in between them. The following function is then defined and minimized for the remaining non-consecutive pairs:
\begin{equation}  \label{equ:objfun1_v1}
f_1(X) = \displaystyle \sum_{i=1}^{n-2} \sum_{j=i+2}^{n} d(x_i,x_j) .
\end{equation}
Minimizing this function is equivalent to forcing all the Euclidean distances $d(x_i,x_j)$, with $j > i + 1$, to be as small as possible. Theoretically, the smallest distance between two amino acids might be 0 according to this function, but this would an unrealistic situation, in which two atoms completely overlap. To avoid this possibility, a constraint can be used to force every pair of non-consecutive amino acids along the chain to have a distance above a certain threshold. Instead of such a constraint, we introduce the function:
\begin{equation}  \label{equ:objfun2}
f_2(X) = \displaystyle \sum_{i=1}^{n-2} \sum_{j=i+2}^{n} exp_+ (th - d(x_i,x_j)).
\end{equation}
Minimizing this function is equivalent to minimizing the number of couples $(x_i,x_j)$ that violates the constraint
\begin{displaymath}
d(x_i,x_j) > th.
\end{displaymath}
The function $exp_+$ in (\ref{equ:objfun2}) corresponds to a function that returns the exponential value of a non-negative argument, and returns zero otherwise. This means that function (\ref{equ:objfun2}) does not have any effect if the corresponding distance is greater than $th$, whereas it grows exponentially as the distance decreases from $th$ to zero. In this work, the value of the threshold $th$ is set to 4.30\AA. The reason for the choice of this particular threshold value is discussed in Section \ref{sec:sim}.

The third and the last requirement we impose is regarding the $\alpha$-helices. We want to have fragments folded as $\alpha$-helices in the protein conformations. As already observed in previous studies \cite{Ceci,MucherinoPF}, the distance between the $C_{\alpha}$ Carbon atom $x_i$ and the $C_{\alpha}$ Carbon atom $x_{i+2}$ can vary within a small interval if both $x_i$ and $x_{i+2}$ are contained in a helix conformation. Therefore, we force the distances $d(x_i,x_{i+2})$ to have a value as close to $c = 5.50{\rm\AA}$ as possible, which is the mean value for such distances that we observe in helix conformations. Then, the third function we consider is as follows:
\begin{equation}  \label{equ:objfun3}
f_3(X) = \displaystyle \sum_{i=1}^{n-2} (d(x_i,x_{i+2}) - c)^2 . \\
\end{equation}
Since the couples $(x_i,x_{i+2})$ are also taken into account in the function $f_1(X)$ (see (\ref{equ:objfun1_v1})), we modify the function $f_1(X)$ as follows:
\begin{equation}  \label{equ:objfun1_v2}
f_1(X) = \displaystyle \sum_{i=1}^{n-3} \sum_{j=i+3}^{n} d(x_i,x_j) .
\end{equation}

Finally, the global optimization problem we solve for simulating protein conformations is as follows:
\begin{equation}  \label{equ:optprobl}
\min_X f(X), \quad f(X) = \gamma_1 f_1(X) + \gamma_2 f_2(X) + \gamma_3 f_3(X) .
\end{equation}
The terms of the objective function are the functions (\ref{equ:objfun1_v2}), (\ref{equ:objfun2}) and (\ref{equ:objfun3}), respectively. The positive weights for each term in the final objective function are given as $\gamma_1$, $\gamma_2$ and $\gamma_3$, respectively. These weights can be altered to change the relative importance of the individual terms. The terms $f_1(X)$ and $f_2(X)$ have a complexity of $O(n^2)$, whereas $f_3(X)$ has a complexity of only $O(n)$, Therefore, $f(x)$ has a complexity of $O(n^2)$, which improves upon the objective function in \cite{Ceci,MucherinoPF}, where the objective function has complexity $O(n^3)$.

%%%%%%%%%%%%%%%%%%%%%%%%%%%%%%%%%%%%%%%%%%%%%%%%%%%%%%%%%%%%%%%%%%%%%
\section{Comparison of the proposed geometric model}  \label{sec:sim}
In this paper, we mainly refer to the studies presented in \cite{Ceci,MucherinoPF} since they involve protein conformations that are simulated using models based on the geometric properties of the proteins. The $\alpha$-helices are simulated by exploiting the concept of protein thickness $\Delta$. It has been proved that, in real $\alpha$-helical conformations, the thickness of a protein can range only within a small interval $I(\Delta) = [2.50,2.80]$. Then, $\alpha$-helices can be simulated by maximizing the number of triplets of amino acids whose thickness fall into this small interval \cite{Ceci}. 

The thickness of a triplet of amino acids can be defined as the radius $r$ of the only circle passing through them, where each amino acid is represented by the spatial coordinates of its $C_{\alpha}$ Carbon atom. In general, the thickness of a whole conformation $X$ containing more than three amino acids is
\begin{displaymath}
\Delta(X) = \min_{i \ne j \ne k \ne i} r(x_i,x_j,x_k),
\end{displaymath}
where $r(x_i,x_j,x_k)$ is the radius of the circles passing through the generic triplet of amino acids of the conformation \cite{Gonzalez}.

We do not consider any requirement based on the protein thickness. However, protein thickness is implicitly considered, and, as we show in Section \ref{sec:experim}, we are still able to simulate protein conformations containing $\alpha$-helices. We prove that when the objective function of the optimization problem (\ref{equ:optprobl}) is minimized, the number of triplets of amino acids having a thickness in the interval $I(\Delta) = [2.50,2.80]$ attains the maximum possible value. 

Minimizing the function $f_1(X)$ is equivalent to minimizing the distances between amino acids, but such distances cannot be below a certain threshold $th$, as controlled by the function $f_2(X)$. Hence, couples of amino acids tend to have a relative distance close to $th$. The following formally proves that, if the couples of amino acids in the triplet $(x_i,x_j,x_k)$ approaches to the threshold $th$, then the thickness of the triplet approaches to $th / \sqrt{3}$. Therefore, the value of the threshold $th$ can be adjusted so that the corresponding thickness value falls into the interval $I(\Delta)$.

\begin{figure}  
\begin{center}
\includegraphics[scale=0.30]{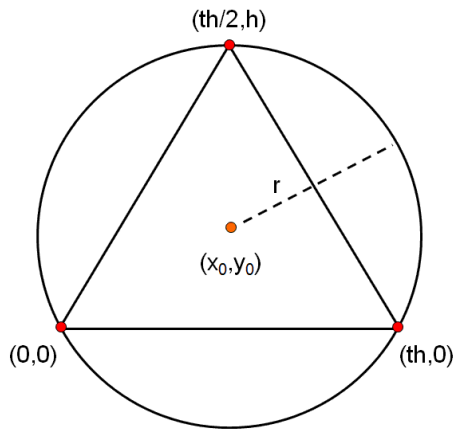}
\end{center}
\caption{A circle passing through the vertices of an equilateral triangle having side $th$ has radius $r = th / \sqrt{3}$.}
\label{fig:triangle}
\end{figure}

Suppose that all possible couples of amino acids in $(x_i,x_j,x_k)$ have a relative distance equal to $th$. As shown in Figure \ref{fig:triangle}, the three amino acids can be seen as the vertices of an equilateral triangle.
For simplicity, the figure is drawn in a two-dimensional space, without losing generality, since three points define a unique plane. One of the amino acids can be put on the origin of the Cartesian system and then the others can have coordinates as shown in Figure \ref{fig:triangle}. The coordinates of these amino acids can be identified, and the passage of a circle among these amino acids can be imposed. The unique circle passing through these amino acids has the center
\begin{displaymath}
(x_0,y_0) = \left( \frac{th}{2},\frac{th}{2\sqrt{3}} \right),
\end{displaymath}
and a radius of 
\begin{displaymath}
r = \frac{th}{\sqrt{3}}.
\end{displaymath}
Thus, if the couples in $(x_i,x_j,x_k)$ approaches a relative distance equal to $th$, then the thickness of the triplet approaches to $th / \sqrt{3}$. It follows that the threshold $th$ must be in the interval $[4.33,4.88]$ for letting the thickness approach to the typical interval $I(\Delta) = [2.50,2.80]$. This proves that the protein thickness is implicitly taken into account in our simplified model. Moreover, this is the reason we consider the value of the threshold as $th = 4.30{\rm\AA}$ in our computational experiments.

%%%%%%%%%%%%%%%%%%%%%%%%%%%%%%%%%%%%%%%%%%%%%%%%%%%%%%%%%%%%%%%%%%%
\section{The Monkey Search algorithm}  \label{sec:ms}

The Monkey Search (MS) is a meta-heuristic approach for global optimization \cite{MucherinoMS,Seref}. It resembles the behavior of a monkey climbing trees in its search for food. The main assumption in this approach is that a monkey is able to survive in a jungle of trees because it is able to remember food sources previously discovered. When the monkey climbs up a new tree for the first time, it can only choose the branches of the tree in a random way, because it does not have any previous experience on that tree. Upon climbing down the tree, the monkey marks tree branches with respect to the quality of the food available in the subtree starting at that branch. When the monkey climbs up the tree again later, using the previous marks on the branches, it tends to choose those branches that lead to the parts of the tree with better quality of food.

In MS, food is represented by the feasible solutions of a global optimization problem. These solutions are stored at the tips of branches of virtual trees that a virtual monkey climbs. The quality of these solutions is evaluated through the value of the objective function of the global optimization problem. A branch also represents a perturbation that creates the solution at the tip of the branch if it is applied to the solution stored at the tip of its parent branch. The trees are not pre-existing, but they are rather built as the monkey climbs up and down. There are three states that the monkey would be in: exploring, climbing down, climbing up.

\emph{Exploring:} At the very beginning, the monkey is located at the root of the tree with a given feasible starting solution. At this point, the tree actually does not exist yet, and therefore its branches need to be created. At each step, two new branches are randomly built. A different perturbation is applied to the current solution, and two new solutions are generated and placed at the tips of the two new branches. At this point, the monkey has to choose between the two branches the one it wants to climb. A random mechanism is used, which is guided by the objective function values of the solutions at the tips of the two branches. Once the chosen branch is climbed, the monkey is located at the tip of that branch. As in the previous case, another two branches are created randomly and the monkey climbs one of them. This procedure continues until the monkey reaches a pre-determined height, which is assumed to be the top of the tree. 

\emph{Climbing Down:} After the monkey reaches the top of the tree, it starts climbing down. The objective function value at the top of the tree is set as the initial marker value. Every time the monkey is at the tip of a branch, first it compares the objective function value of that branch to the current marker value. If the objective function value is better than the current marker value, then the marker is updated with the objective function value of the branch; otherwise the marker value is kept as is. Then, the branch is marked with the marker value and the monkey climbs down to the parent branch. Note that the marker value on a branch is equivalent to the best possible value in the subtree that starts from that branch and expands up. These marker values are later used to guide the monkey on its way up to the unexplored parts of the tree that potentially have better food resources. While the monkey is climbing down, it randomly chooses to climb up again without necessarily reaching the root of the tree.

\emph{Climbing Up:} When the monkey decides to restart climbing up, it first uses the existing branches until it reaches a branch that is at the border between the explored and the unexplored part of the tree. While climbing up the existing branches, whenever the monkey is at the tip of a branch, it chooses between the two children branches according to a random process based on the previous marks left while the monkey was in the climbing down mode. Naturally, the probabilistic mechanism is in favor of choosing a branch with a better mark. Since this mechanism is based on probabilities, the monkey still has a chance to follow the worse branch. Whenever the monkey reaches the undiscovered branches it switches to \emph{exploring} mode and continues exploring the undiscovered branches until it reaches the top of the tree.

MS is based on a set of parameters that influence the convergence of the algorithm. The \emph{height} of the trees is the total number of branches that the monkey can climb from the root to the top. The \emph{number of climbing up} allowed is equivalent to the number of times the monkey decides to climb up before it stops searching a tree. \emph{Memory size} determines number of best solutions found on each individual tree in order to avoid local minima. Each time the monkey reaches the maximum number of paths allowed, it starts climbing a new tree starting with a different root solution, which is usually set as a combination of the solutions in the memory. \emph{Number of starting trees} sets the number of initial trees with completely random starting points. The motivation behind multiple starting trees is to explore different parts of the global optimizaiton problem's domain. The memory is updated with the best solutions found from these starting trees directly before any consequent tree is formed using the solutions in the memory. MS stops when all of the solutions in the memory of that tree are sufficiently close to one another. 

The perturbations applied when the monkey is in exploring mode play an important role. They can be totally random, as it is done in standard Simulated Annealing algorithms \cite{Kirkpatrick}, but they can also be customized and tailored to the global optimization problem. The idea is to exploit strategies borrowed from other meta-heuristic search methods, such as Genetic Algorithms \cite{Goldberg}, the Harmony Search \cite{Geem}, etc., and/or develop perturbations which are inspired by the problem itself. If a certain set of perturbations is defined, then the perturbation to be applied for generating a new branch can be chosen with a uniform probability or it can adaptively be redirected towards the succesful ones.

%%%%%%%%%%%%%%%%%%%%%%%%%%%%%%%%%%%%%%%%%%%%%%%%%%%%%%%%%%%%%%%%%%%%%%%%%%%%%%%%%%%%%
\section{Computational experiments}  \label{sec:experim}

In this section, we present computational experiment performed with the aim of simulating a large set of protein conformations using the model discussed in Section \ref{sec:opt}. The optimization problem (\ref{equ:optprobl}) is solved using Monkey Search, which is described in Section \ref{sec:ms}. We point out that our protein conformations are represented in two different ways. Since the three terms in the objective function $f(X)$ of the optimization problem (\ref{equ:optprobl}) are based on the Euclidean distances between couples of $C_{\alpha}$ Carbon atoms, it is natural to use
\begin{displaymath}
X = \{ x_1,x_2,\dots,x_n \} ,
\end{displaymath}
to represent the proteins, where each $x_i$ is the three-dimensional vector of coordinates associated with the position of a $C_{\alpha}$. However, this representation $X$ is not so efficient since it provides a conformation with more arrangement possibilities than it actually can. For instance, if no constraints are used, then $x_i = x_j$ might be possible for some $i \ne j$, and this is obviously unrealistic. Although it is easy to evaluate $f(X)$ with this representation, we prefer the following representation 
\begin{displaymath}
\bar{X} = \{ \Phi_3,\Phi_4,\dots,\Phi_n;\Psi_3,\Psi_4,\dots,\Psi_n \}, 
\end{displaymath}
which is based on the so-called {\it dihedral angles}. Along the backbone of a protein conformation, three dihedral angles can be defined, usually denoted by $\Phi$, $\Psi$ and $\omega$, and one of them, $\omega$, is almost always constant and equal to $\pi$. It is well-known that the representation $\bar{X}$ suffices for building the three-dimensional coordinates of the whole protein backbone, including the $C_{\alpha}$ Carbon atoms. Moreover, unrealistic conformations can easily be avoided using this representation, because there may be couples of angles ($\Phi$,$\Psi$) that have never been observed in nature. Therefore, we use the representation $\bar{X}$ for the optimization process, where the generic $\Phi_i$ and $\Psi_i$ are bounded to have naturally occuring values, and the objective function is evaluated after a transformation from $\bar{X}$ to $X$.

The Monkey Search algorithm is used for solving the problem (\ref{equ:optprobl}) and simulating 1000 protein conformations. The software procedures for MS are written in C. The set of parameters for the Monkey Search are tuned by performing preliminary computational experiments. The \emph{height} of a tree is set to 40 branches, and the \emph{number of climb-ups} allowed is set to 20. The \emph{memory size} is set to 10 to keep the best 10 conformations that are used to generate starting solutions for the consequent trees. The \emph{number of starting trees} is set to 100, for the first 100 trees that use random conformations as their starting solutions. The stopping criteria is usually reached after 3000 trees are explored. Each simulation takes approximately 20 minutes on a 2GHz Intel processor, running Windows XP.

All the simulated conformations are formed by 65 three-dimensional points $x_i$ representing the position of a $C_{\alpha}$ Carbon atom, and hence the position of an amino acid. We chose 65 as length of our conformations because some of the smallest all-$\alpha$ proteins currently known have 65 amino acids. We chose 8 proteins as targets in our comparisons. The labels of such proteins in the Protein Data Bank (PDB) \cite{Berman,PDB} are {\tt 1ame}, {\tt 1b7i}, {\tt 1ekl}, {\tt 1gzi}, {\tt 1jc6}, {\tt 1v66}, {\tt 2alc} and {\tt 2cro}. All these proteins are formed by exactly 65 amino acids and they are all-$\alpha$ proteins. The weights $\gamma_1$, $\gamma_2$ and $\gamma_3$ of the function $f(X)$ are also tuned by preliminary computational experimets. Since the triplets of values reported in Table \ref{tab:weights} produces better results compared to other combinations, we use these values for generating a large set of conformations. 
\begin{table}  
\begin{center}
\begin{tabular}{|c|c|c|c|}
\hline
{} & $\gamma_1$ & $\gamma_2$ & $\gamma_3$ \\
\hline\hline
test1 & 0.1 & 1.2 & 0.1 \\
\hline
test2 & 0.1 & 1.0 & 0.1 \\
\hline
test3 & 0.5 & 1.0 & 0.1 \\
\hline
test4 & 0.1 & 1.0 & 0.5 \\
\hline
\end{tabular}
\end{center}
\caption{The four triplets of weights used in the computational experiments.}
\label{tab:weights}
\end{table}
The 1000 conformations we simulated can be divided into 4 groups of 250, with each group corresponding to one of the triplets of the weights $\gamma$ reported in Table \ref{tab:weights}. Since the Monkey Search is a meta-heuristic algorithm, many different solutions to the optimization problem can be found by changing the seed of the random number generator. The number of local minima that the objective function $f(X)$ may have is unknown, and the search may actually stop at one of these local minima. However, we believe that, since the proposed model is quite simplified, the objective function may contain several global optima, and the Monkey Search cannot discriminate among them. All of the found solutions satisfy the imposed requirements. 

Our computational results show that the majority of the simulated conformations have no similarities with the target proteins. This was expected, because the proposed model is very simple and does not consider information regarding the amino acid sequence of the proteins. However, there are simulated conformations that are similar to the target proteins. The similarities between the conformations are evaluated by using the software {\tt ProFit} \cite{Martin}, which computes the Root Mean Square Deviation (RMSD) between the target and the simulated conformations. Only the backbone atoms of the target proteins are considered for the comparison so that the compared conformations have the same number of atoms. 

In general, the RMSD values between protein conformations having the same length and the same secondary structures is large. However, there are conformations we simulated which correspond to relatively small RMSD values if compared to one of the target proteins. For instance, the RMSD value between one of the simulated conformations and the {\tt 2cro} protein is 5.36\AA. Another simulated conformation correspond to an RMSD value of 4.97\AA\ if it is compared to the protein {\tt 1v66}. Even though these are just two examples, in which the RMSD values are not large (the mean RMSD value obtained is about 10\AA), we showed that our model is able to find conformations that are {\it similar} to the ones that actually exist in nature. This result is quite surprising. Indeed, the model we use is very simple, if it is compared to the models used for performing protein predictions in general. In fact, our model is very far from predicting a protein conformation. What is interesting is that our model can generate, in a reasonable amount of time, a very large set of protein conformations in which {\it there usually is} a conformation that is similar to one a natural protein may have. The larger is the number of simulated conformations, the higher is the probability to generate a conformation very similar to a target protein.

\begin{figure}  
\begin{center}
\begin{tabular}{|c|c|c|}
\hline
\includegraphics[scale=0.28]{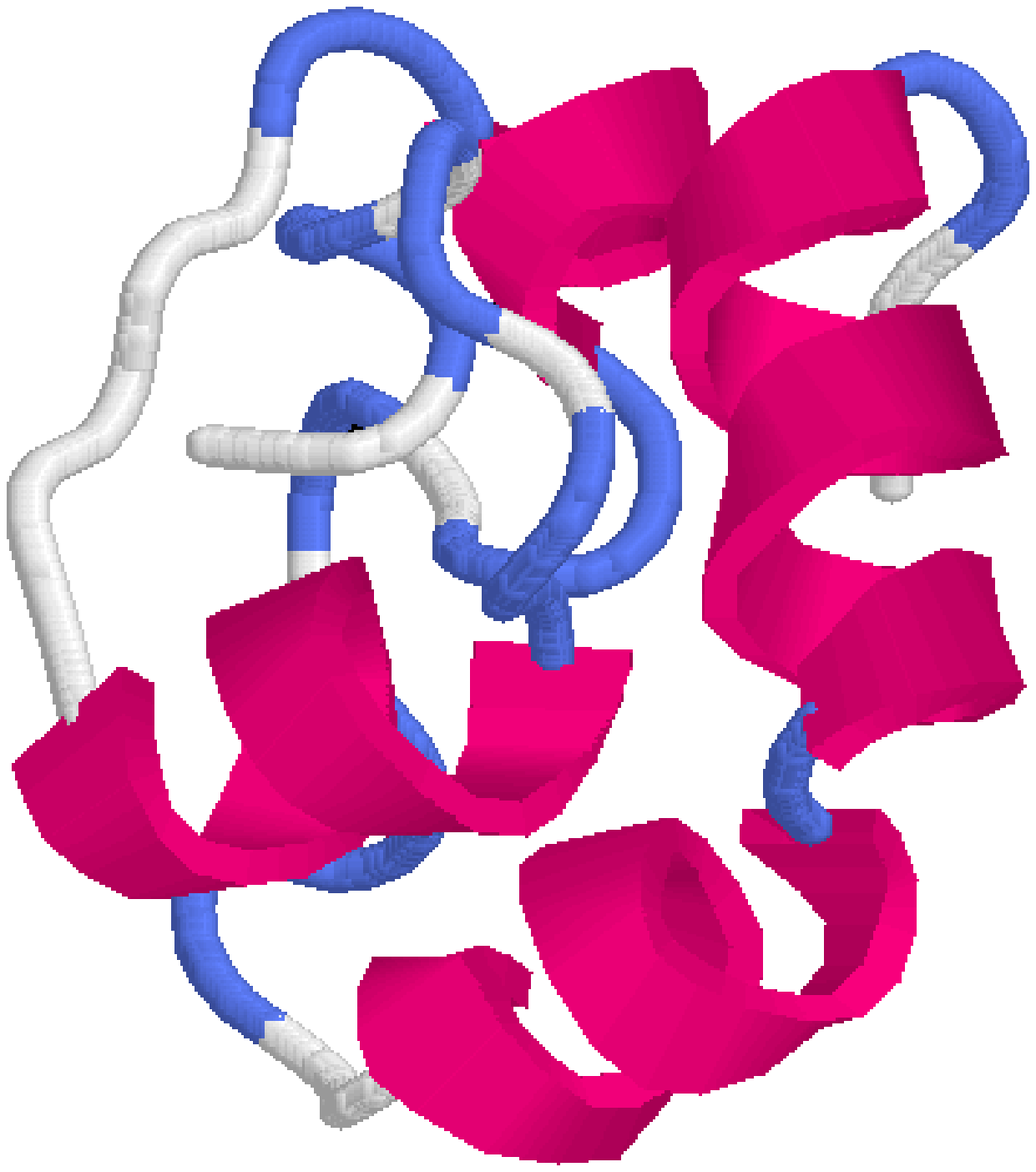} &
\includegraphics[scale=0.28]{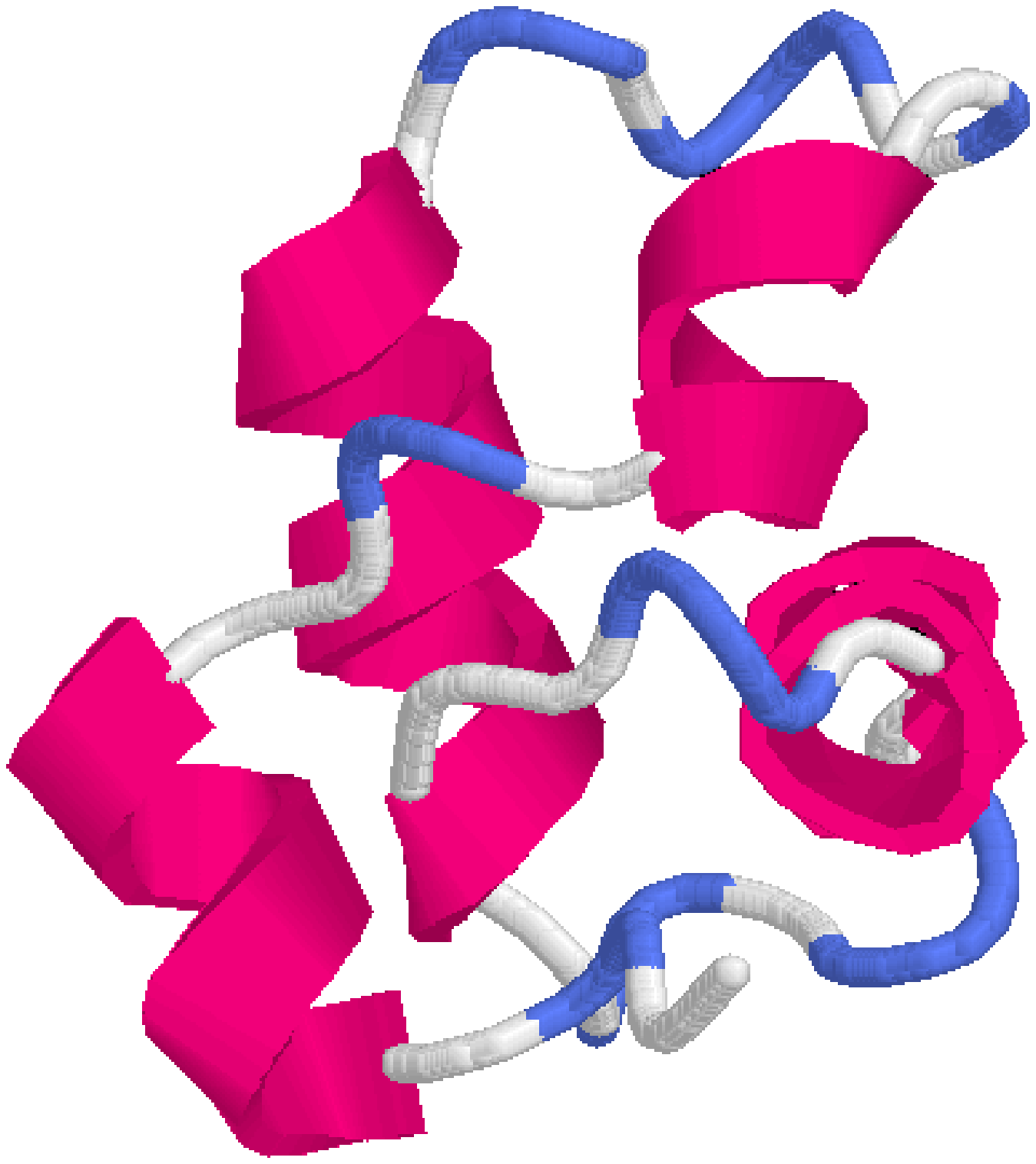} &
\includegraphics[scale=0.28]{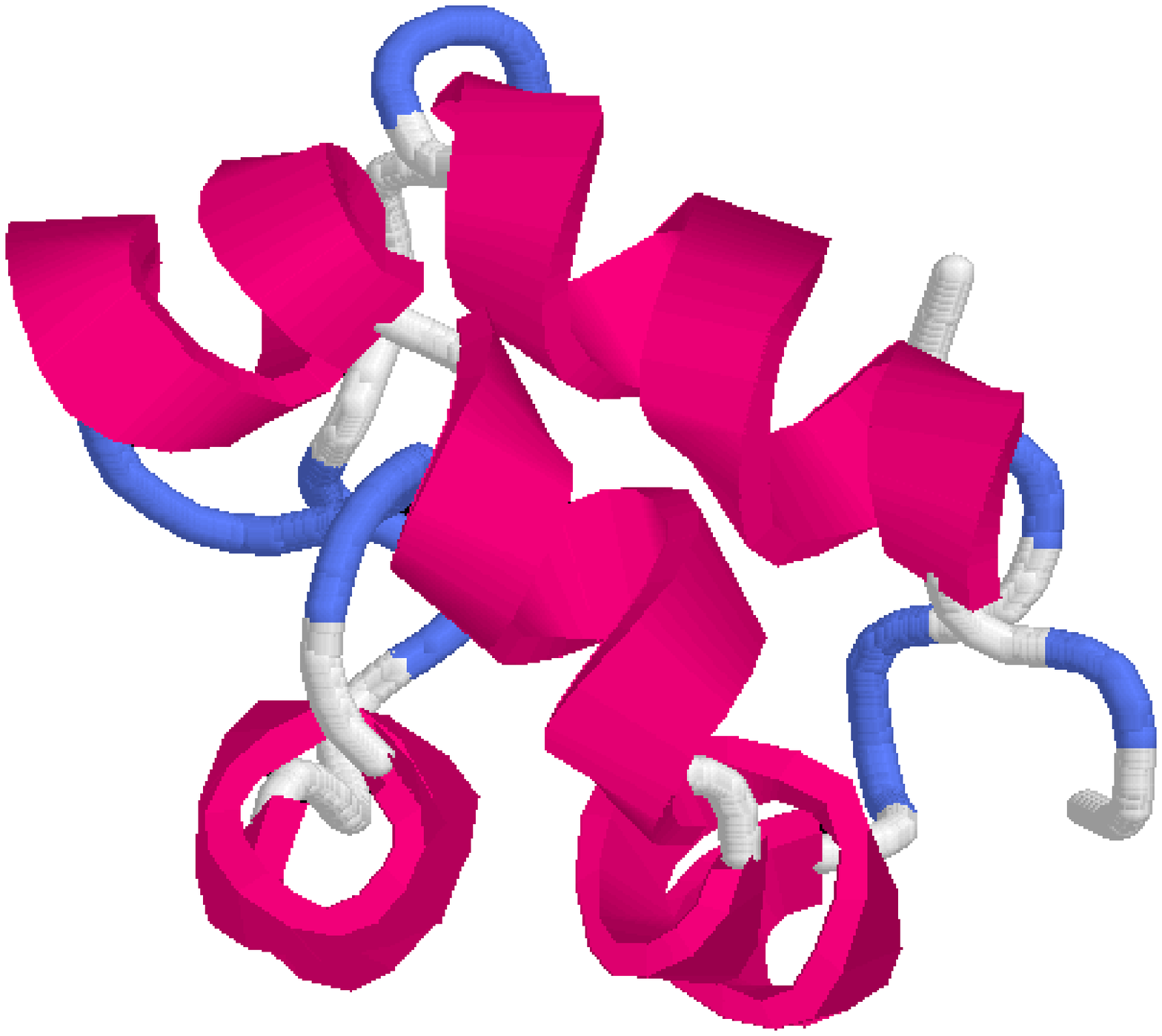} \\
{\tt (a)} & {\tt (b)} & {\tt (c)} \\
\hline
\end{tabular}
\end{center}
\caption{Three conformations simulated by the considered model.}
\label{fig:simpro}
\end{figure}

In Figure \ref{fig:simpro}, three simulated conformations are shown.
The conformation in Figure \ref{fig:simpro}(a) is the most similar to the protein {\tt 2cro}. The corresponding RMSD value is 5.36\AA\ and the triplet of weights used for generating it is $(\gamma_1,\gamma_2,\gamma_3) = (0.1,1.0,0.1)$. The conformations in Figure \ref{fig:simpro}(b) and Figure \ref{fig:simpro}(c) are more similar to the protein {\tt 1v66}. In particular, the first one and {\tt 1v66} have an RMSD value of 5.16\AA, whereas the second one and {\tt 1v66} have an RMSD value of 4.97\AA. The triplet of weights $(\gamma_1,\gamma_2,\gamma_3)$ is $(0.1,1.0,0.5)$ for the first one and $(0.1,1.0,0.1)$ for the second one. Note that these conformations satisfy the imposed geometric requirements, and therefore they are solutions of the formulated optimization problem. They are indeed compact, their backbones never intersect themselves, and helices are formed. Even though all the simulated conformations satisfy these requirements, only the ones in Figure \ref{fig:simpro} are similar to real protein conformations.

%%%%%%%%%%%%%%%%%%%%%%%%%%%%%%%%%%%%%%%%%%%%%%%%%%%%%%%%%%%%%%%%%%%%%%%%%%%%%%%%
\section{Conclusion}  \label{sec:concl}

A model for simulating protein conformations is discussed in this paper. The basic idea is to impose simple geometric requirements on a virtual chain of amino acids. A global optimization problem is formulated whose solutions provide conformations meeting these geometric requirements. We force our simulated conformations to be compact and to contain fragments folded in $\alpha$-helices. This model is inspired from the works recently published in \cite{Banavar02,Ceci,MucherinoPF}, where the so-called {\it tube model} and its variances have been proposed. We show that the protein thickness, which is the main parameter involved in these previous models, is implicitly considered in our simple model.

Computational experiments show that conformations having the typical properties of the all-$\alpha$ proteins can be simulated using the proposed model. Monkey Search algorithm is used to solve the global optimization problem based on the proposed model and allowed the simulation of a high-quality set of conformations in a reasonable amount of time. Some of these simulated conformations are similar to real proteins chosen as targets, measured by the the root mean square deviation (RMDS) between the simulated proteins and the actual proteins. This can be considered as a significant result, because the conformations are simulated by a very simple model that considers few geometric requirements only.

The results discussed in this paper suggest that the geometric properties of the proteins can be exploited in order to simulate, or even predict protein conformations. The model discussed in this paper takes few geometric requirements into account, and the introduction of other requirements may improve the quality of the generated conformations. The more the requirements are imposed, the more the set of simulated conformations should converge toward a subset, in which only conformations close to the real ones are contained. Naturally, our model could never be used for protein prediction if information regarding the chemical composition of the amino acid sequence is not taken into consideration. However, we think that geometric properties provide a valid alternative to the force fields that are currently used in prediction models. Since standard physics-based models for protein prediction are still far from finding a good and accurate protein conformations, we believe it may be worth to follow this uncommon approach.

\end{document}